\magnification=\magstep1
\input amstex
\documentstyle{amsppt}

\TagsOnRight
\rightheadtext{A Trinomial Analogue of Bailey's Lemma}
\leftheadtext{G. E. Andrews and A. Berkovich}
\NoBlackBoxes

\def\fr{\frac}
\def\D{\Cal D}
\def\wp{\widehat{p}\,}
\def\wr{\widehat{r}}
\def\ws{\widehat{s}}

\def\wchi{\widehat{\chi}}
\def\wtr{\widetilde{r}}
\def\wts{\widetilde{s}}
\def\wtp{\widetilde{p}\,}
\def\wtin{\widetilde{\in}}
\def\wtchi{\widetilde{\chi}\,}
\def\wbeta{\widetilde{\beta}}
\def\walpha{\widetilde{\alpha}}
\def\wdelta{\widetilde{\delta}}
\def\wgamma{\widetilde{\gamma}}
\def\sumlo{\sum_{L=0}^{\infty}}
\def\sumro{\sum_{r=0}^{\infty}}

\topmatter
\title 
	A Trinomial Analogue of Bailey's Lemma and $\bold{N = 2}$
	Superconformal Invariance
\endtitle
\author
	G. E. Andrews$^{1)}$  \\
	A. Berkovich$^{2)}$ \\ \\ \\
       \it{Dedicated to Dora Bitman on her 70th birthday} \\
\endauthor
\abstract
	We propose and prove a trinomial version of the celebrated
	Bailey's lemma.  As an application we obtain new fermionic 
	representations for characters of some unitary as well as
	nonunitary models of $N = 2$ superconformal field theory
	(SCFT).  We also establish interesting relations between
	$N = 1$ and $N = 2$ models of SCFT with central charges
	$\fr32 \left( 1 - \fr{2(2 - 4\nu)^2}{2(4\nu)}\right)$ and
	$3\left(1 - \fr2{4\nu}\right)$. A number of new mock theta 
        function identities are derived.
\endabstract
\endtopmatter
\document

\footnote"{}"{$^{1)}$e-mail: \ andrews\@math.psu.edu}
\footnote"{}"{$^{2)}$e-mail: \ berkov\_a\@math.psu.edu}
\footnote"{}"{Partially supported by National Science Foundation
Grant: \ DMS-9501101.}

\baselineskip 20pt

\subhead
	1. \ Brief review of Bailey's method and its generalizations
\endsubhead

It may come as a surprise that Manchester, England was an ideal setting
for pure mathematics during the height of World War II.  However, a 
variety of historical coincidences conspired to make this the case.
In particular, mathematics that would later prove extremely valuable
in the development of statistical mechanics and conformal field
theory (CFT) flourished there.

Essentially, Bailey, extending the original ideas of Rogers, came up
with a new method $[1,2]$ of deriving Rogers-Ramanujan type 
identities during the winter 1943--44.  Hardy who was then editor
for the Journal of London Mathematical Society sent a referee
report with Dyson's name on it back to Bailey.  Bailey's reply
was immediate.  A charming account of Dyson-Bailey collaboration
appears in Dyson's article, A Walk Through Ramanujan's Garden [3]. 

A few years later, Slater, in a study building on Bailey's work,
systematically derived 130 identities of Rogers-Ramanujan type [4,5].
In the last decade, Bailey's technique was streamlined and 
generalized by Andrews [6] and further extended by Agarwal,
Andrews and Bressoud [7,8].

Bailey's method may be summarized as follows.  Let $\alpha = \{\alpha_r
\}_{r \geq 0}$, $\beta = \{\beta_L \}_{L\geq 0}$ be sequences related
by identities
$$
	\beta_L = \sum_{r=0}^{\infty} \fr{\alpha_r}{(q)_{L-r}(aq)_{L+r}}
	\quad,\quad L \in Z_{\geq 0}
\tag1.1
$$
$$
	 (a)_n = \cases
		 \infty \,, \qquad & n \in Z_{<0}    \\	
		 1\,, \qquad & n = 0   \\
		 (1-a)(1 - aq) \cdots (1 - aq^{n-1})\,, \qquad & n \in Z_{>0}
	 \endcases
 \tag1.2
 $$
 and $\gamma = \{\gamma_L\}_{L\geq 0}$, $\delta = \{ \delta_{r}\}_{r\geq 0}$
 be another pair of sequences related by 
 $$
	 \gamma_L = \sum_{r=L}^{\infty} 
		 \fr{\delta_r}{(q)_{r-L}(aq)_{r+L}}.
 \tag1.3
 $$
 Then the new identity
 $$
	 \sum_{L=0}^{\infty} \alpha_L \gamma_L =
		 \sum_{L=0}^{\infty} \beta_L \delta_L
 \tag1.4
 $$
 holds.  A pair of sequences $(\alpha,\beta)$ that satisfies (1.1)
 is called a Bailey pair relative to $a$.  Analogously, a pair of
 sequences $(\gamma,\delta)$ subject to (1.3) is referred to as 
 conjugate Bailey pair relative to $a$.  In [2], Bailey proved that
 $$
 \align
	 \gamma_L & = \fr{(\rho_1,\rho_2)_L (aq/\rho_1 \rho_2)^L}
		 {(aq/\rho_1,aq/\rho_2)_L} \;\;\fr1{(q)_{M-L}(aq)_{M+L}}
		 \tag1.5
	 \\
	 \delta_L & = \fr{(\rho_1,\rho_2)_L (aq/\rho_1 \rho_2)^L}
		 {(aq/\rho_1,aq/\rho_2)_M} \;\;\fr{(aq/\rho_1 \rho_2)_{M-L}}
		 {(q)_{M-L}}
		 \tag1.6
 \endalign
 $$
 with $(a_1,a_2)_L \equiv(a_1)_L (a_2)_L$, $L \leq M \in Z_{\geq 0}$
 satisfy (1.3) for any choice of parameters $\rho_1,\rho_2$.  Combining
 (1.4, 1.5, 1.6) yields
 $$
 \multline
	 \sum_{L=0}^{\infty} \fr{(\rho_1,\rho_2)_L}{(aq/\rho_1,aq/\rho_2)_M}
	 \left( \fr{aq}{\rho_1 \rho_2}\right)^L 
	 \fr{(aq/\rho_1 \rho_2)_{M-L}}{(q)_{M-L}} \beta_L =   \\ 
	 \sum_{L=0}^{\infty} \fr{(\rho_1,\rho_2)_L(aq/\rho_1 \rho_2)^L}
		 {(aq/\rho_1,aq/\rho_2)_L} \quad \fr{\alpha_L}{(q)_{M-L}
		 (aq)_{M+L}}.
 \endmultline
 \tag1.7
 $$
 From the last equation, we deduce immediately:

 {\bf (Bailey's lemma)}  Sequences $(\alpha',\beta')$ defined by
 $$
	 \alpha_L'  = \fr{(\rho_1,\rho_2)_L (aq/\rho_1 \rho_2)^L}{(aq/\rho_1,
		 aq/\rho_2)_L} \;\alpha_L
 \tag1.8
 $$
 $$
	 \beta_L'  = \sum_{r=0}^{\infty}
		 \fr{(\rho_1,\rho_2)_r (aq/\rho_1 \rho_2)^r(aq/\rho_1 
			 \rho_2)_{L-r}}
		 {(aq/\rho_1,aq/\rho_2)_L(q)_{L-r}} \;\beta_r
 \tag1.9
 $$
 form again a Bailey pair relative to $a$.  Obviously, Bailey's 
 lemma can be iterated and infinitum leading to a ``Bailey chain''
 [6,9] of new identities
 $$
	 (\alpha,\beta) \to (\alpha',\beta') \to (\alpha'',\beta'')
	 \to (\alpha''',\beta''') \to \dots
 \tag1.10
 $$
 with parameter $a$ remaining unchanged throughout the chain.  The
 notion of a Bailey chain was upgraded to a  ``Bailey lattice'' in
 [7,8] where it was shown how to pass from a Bailey pair with given
 parameter $a$ to another pair with arbitrary new parameter.

 Further important developments have taken place in the last few
 years.
 In [10], Milne and Lilly found higher-rank generalizations of
 Bailey's lemma. Many new polynomial identities of Rogers-Ramanujan type
 were discovered in [11-18] as result of recent progress in CFT and
 Statistical Mechanics initiated by the Stony Brook group [19-21].
 Following observation made by Foda-Quano [22], these identities were recognized
 as new $(\alpha, \beta )$ pairs. New $(\gamma,\delta )$ pairs were discovered
 in [23,24]. Intriguing connections between Bailey's lemma and the so-called
 renormalization group flows connecting different models at CFT were
 discussed in [17,25-27].

 This paper is intended as the first step towards a multinomial (or
 higher-spin) generalization of Bailey's lemma.  Here we concentrate
 on the trinomial case.  Our main assertion is
 \proclaim
 {Theorem 1}  (Trinomial analogue of Bailey's lemma)

 If for $L \geq 0$, $a = 0,1$
 $$
	 \wbeta_a(L) = \sum_{r=0}^L \walpha_a(r) 
		 \fr{T_a(L,r,q)}{(q)_L}
 \tag1.11
 $$
 then for $M \in Z_{\geq 0}$
 $$
	 \sumlo (-1)_L q^{\fr{L}{2}} \wbeta_0 (L) = \sumro\walpha_0
	 (r) \fr{(-1)_{M+1}}{q^{\fr{r}{2}} + q^{-\fr{r}{2}}} \quad
	 \fr{T_1(M,r,q)}{(q)_M}
 \tag1.12
 $$
 and
 $$
 \aligned
	 \sumlo  (-q^{-1})_L q^L \wbeta_1 (L) =
	 \sumro  \walpha_1(r) & \fr{(-1)_M}{(q)_M}\bigg\{ T_1(M,r,q)  \\
	 & - \fr{(1 - q^M)}{1 + q^{-1-r}} T_1(M-1,r+1,q)    \\
	 & \left. - \fr{(1 - q^M)}{1 + q^{r-1}} T_1(M-1,r-1,q)\right\}
 \endaligned
 \tag1.13
 $$
 where $T_a(L,r,q)$ are $q$-trinomial coefficients [28] to be defined
 in the next section.  The pair of sequences $(\walpha_a,\wbeta_a)$
 that satisfies identities (1.11) will be called a trinomial Bailey
 pair.
 \endproclaim

 The rest of this paper is organized as follows.  In Section 2 we shall
 collect the necessary background on $q$-trinomials and then prove
 Theorem 1.  In Section 3 we shall exploit this theorem to derive a
 number of new $q$-series identities related to characters of $N = 2$
 SCFT. We conclude with a brief
 discussion of the physical significance of our results and some
 comments about possible generalizations.

 \subhead
	 2. \ $\bold{q}$-Trinomial Coefficients and a Trinomial Analogue 
	 of Bailey's Lemma
 \endsubhead

 \subhead
 2.1 \ Preliminaries
 \endsubhead

 Before turning our attention to the $q$-trinomial coefficients, 
 let us briefly recall the ordinary trinomials $\binom{L}{A}_2$
 [28], defined by 
 $$
	 \left( x + 1 + \fr1{x}\right)^L = \sum_{A = -L}^L
	 \binom{L}{A}_2 x^A
 \tag2.1
 $$
 $$
	 \binom{L}{A}_2 = 0\,, \qquad |A| > L.
 \tag2.2
 $$
 By applying the binomial theorem twice to (2.1) we find by 
 coefficient comparison that
 $$
	 \binom{L}{A}_2 = \sum_{j \geq 0} 
	 \fr{L!}{j!(j + A)! \;(L - 2j - A)!}.
 \tag2.3
 $$
 Furthermore, it is easy to deduce from (2.1) the following recurrences
 $$
	 \binom{L}{A}_2 = \binom{L - 1}{A - 1}_2 +
	 \binom{L - 1}{A}_2 + \binom{L - 1}{A + 1}_2
 \tag2.4
 $$
 which along with (2.2) and
 $$
	 \binom{0}{0}_2 = 1
 \tag2.5
 $$
 specify trinomials uniquely.  Equations (2.4) and (2.5) lead to the
 Pascal-like triangle for $\binom{L }{A}_2$ numbers
 $$
 \matrix
	 & & & & &  1 & &    \\
	 & & & &  1 & 1 & 1   \\ 
	 & & & 1 & 2 & 3 & 2 & 1    \\
	 & & 1 & 3 & 6 & 7 & 6 & 3 & 1  \\
	 & -  & - & - & - & - & - & -  & - & -   \\
	 - & -  & - & - & - & - & - & - & - & - & -\;\;.
 \endmatrix 
 $$

 $q$-Analogues of trinomials were forced into existence in the work 
 of Andrews and Baxter on the generalized Hard Hexagon model [28].
 These analogues proved to play an important role in Partition
 Theory [29-31] and Statistical mechanics [32-34].  
 Unlike binomials, trinomials admit not one but many $q$-analogues
 which we now proceed to describe.

 \subhead
 2.2 \ Definitions and Properties of $\bold{q}$-trinomials
 \endsubhead 

 The straightforward $q$-deformation of (2.3) is as follows
 $$
	 \binom{L;B;q}{A}_2 = \sum_{j \geq 0} 
	 \fr{q^{j(j+B)} (q)_L}{(q)_j (q)_{j+A} (q)_{L-2j-A}}.
 \tag2.6
 $$
 Let us further define another useful $q$-analogue of (2.3):
 $$
	 T_n(L,A,q) = q^{\fr{L(L-n) - A(A-n)}{2}}
	 \binom{L;A-n;q^{-1}}{A}_2\,, \quad n \in Z.
 \tag2.7
 $$ 
 Polynomials $T_n(L,A,q)$ are symmetric under $A \to - A$
 $$
	 T_n(L,A,q) = T_n (L,-A,q)
 \tag2.8
 $$
 and vanish for $|A| > L$:
 $$
	 T_n(L,A,q) = 0 \quad \text{ if } \quad |A| > L\,.
 \tag2.9
 $$
 The generalization of Pascal-triangle type recurrences (2.4)
 found in [31, 33] is
 $$
 \align
	 T_n(L,A,q) & = T_n(L-1,A-1,q) + T_n(L - 1,A+1,q) + \\
	 & \qquad q^{L - \fr{1 + n}{2}} T_n (L - 1,A,q) + 
		 (q^{L-1}-1)T_n(L-2,A,q)\,.
 \tag2.10
 \endalign
 $$
 Additionally, there are four more identities needed
 $$
 \align
	 T_n(L,A,q) & = T_{n+2}(L,A,q) + (q^L - 1)q^{-\fr{1 + n}{2}}
	 T_n(L - 1,A,q)   \tag2.11    \\
	 q^{\fr{L-A}{2}} T_{n+1}(L,A,q) & = T_n(L,A,q) + (q^L - 1)
		 T_n(L - 1,A + 1,q)   \tag2.12  \\
	 T_1(L,A,q) & - T_1(L - 1,A,q) =   \tag 2.13    \\
		 & = q^{\fr{L + A}{2}} T_0 (L-1,A+1,q) + q^{\fr{L-A}{2}}
		 T_0(L-1,A-1,q)   \\
	 T_1(L,A,q) & + T_1(L - 1,A,q) =   \tag 2.14    \\
	  & = T_{-1}(L - 1,A + 1,q) + T_{-1}(L - 1,A- 1,q) + 2T_{-1}
		 (L - 1,A,q).
 \endalign
 $$
 Identities (2.11), (2.12), (2.13) follow from equations (2.24), 
 (2.23), (2.16) of [28] and identity (2.14) is equation (4.5)
 of [33].

 Next we shall require the limiting formula
 $$
	 \lim_{L \rightarrow \infty} T_1(L,A,q) =	
		 \fr{(-q)_{\infty}}{(q)_{\infty}}
 \tag2.15
 $$
 which is equation (2.51) of [28].

 Let us combine (2.11) and (2.12) with $n = -1$ to obtain 
 $$
 \align
	 q^{\fr{L-A}{2}} T_0 (L,A,q)  & - T_1(L,A,q) =  \tag2.16  \\
	 & = (q^L - 1)\{T_1(L-1,A + 1, q) + T_{-1}(L-1,A,q)\}.
 \endalign
 $$
 We now replace $A$ by $-A$ in the above equation to get, with the 
 help of (2.8),
 $$
 \align
	 q^{\fr{L+A}{2}} T_0 (L,A,q)  & - T_1(L,A,q) =  \tag2.17  \\
	 & = (q^L - 1)\{T_1(L-1,A - 1, q) + T_{-1}(L-1,A,q)\}.
 \endalign
 $$
 If we add (2.16) and (2.17) and use (2.14), the result is
 $$
 \align
	 q^{\fr{L}{2}} \left( q^{\fr{A}{2}} + q^{-\fr{A}{2}}\right)
	 T_0 (L,A,q)  & - 2 T_1(L,A,q) =  \tag2.18  \\
	 & = (q^L - 1)\{T_1(L,A,q) + T_1(L-1,A,q)\}
 \endalign
 $$
 which may be conveniently rewritten as
 $$
 \align
	 T_0 (L,A,q)  & = \fr{q^{-\fr{L}{2}}}{q^{\fr{A}{2}} + q^{-\fr{A}{2}}}
		 (1 + q^L)  T_1(L,A,q)     \tag2.19  \\
	 & \quad- \fr{q^{-\fr{L}{2}}(1 - q^L)}{q^{\fr{A}{2}} + q^{-\fr{A}{2}}}
		 T_{1}(L-1,A,q).
 \endalign
 $$
 We're now ready to prove Theorem 1.

 \subhead
 2.3 \ Proof of Theorem 1
 \endsubhead

 We shall prove Theorem 1 in two steps.  First, we find a trinomial
 analogue of a conjugate Bailey pair and then use a standard Bailey
 Transform argument.  To this end let us introduce an auxiliary
 function $\phi(L,q)$
 $$
	 \phi(L,q) = q^{\fr{L}{2}} \fr{(-1)_L}{(q)_L}
 \tag2.20
 $$
 with the easily verifiable property
 $$
	 \phi(L + 1,q) = \sqrt{q} \fr{1 + q^L}{1 - q^{L + 1}}
	 \phi(L,q)   
 \tag2.21
 $$
 Next we multiply both sides of (2.19) by $\phi(L,q)$ and sum 
 both extrems of the result on $L$ from $A$ to $M$ to obtain 
 with the aid of (2.21) the following
 $$
	 \sum_{L=A}^{M}  q^{\fr{L}{2}} \fr{(-1)_L}{(q)_L}
	 T_0(L,A,q) = \fr{(-1)_{M + 1}}{(q)_M} \;
	 \fr{T_1(M,A,q)}{q^{\fr{A}{2}} + q^{-\fr{A}{2}}}
 \tag2.22
 $$
 which can be restated as a trinomial analogue of the conjugate
 relation (1.3)
 $$
	 \wgamma_0(A,M) = \sum_{L=A}^{\infty}  \wdelta_0(L,M)\;
	 \fr{T_0(L,A,q)}{(q)_L}
 \tag2.23
 $$
 with conjugate pair $(\wgamma_0,\wdelta_0)$ defined as
 $$
	 \wgamma_0(A,M) = \fr{(-1)_{M+1}}{q^{\fr{A}{2}} + q^{-\fr{A}{2}}}\;
	 \fr{T_1(M,A,q)}{(q)_M}
 \tag2.24
 $$
 $$
	 \wdelta_0(L,M) = \theta(L \leq M) q^{\fr{L}{2}} (-1)_L
 \tag2.25
 $$
 where
 $$
	 \theta(a \leq b) = \cases
		 1 \qquad & \text{ if } a \leq b     \\
		 0 \qquad & \text{ otherwise }.
	 \endcases
 \tag2.26
 $$
 The proof of the first statement of Theorem 1 (1.12) now easily
 follows by a Bailey Transform argument
 $$
 \gathered
	 \sumro \walpha_0(r)\wgamma_0(r,M) = \sumro \walpha_0(r) 
	 \sum_{L=r}^{\infty} \wdelta_0(L,M) \fr{T_0(L,r,q)}{(q)_L} 
	 \\
	 = \sumlo \wdelta_0(L,M) \sum_{r=0}^L \walpha_0(r)
	 \fr{T_0(L,r,q)}{(q)_L} = \sumlo \wdelta_0(L,M) \wbeta_0(L). 
 \endgathered
 \tag2.27
 $$

 Substituting (2.24) and (2.25) into (2.27) we arrive at the
 desired result (1.12).

 Similar to the binomial case, identity (1.12) can be interpreted
 as a defining relation for new trinomial Bailey pair $\left(\walpha_1,
 \wdelta_1\right)$.  However, unlike the binomial case, the second analogue
 of (1.3)
 $$
	 \wgamma_1(A,M) = \sum_{L=A}^{\infty} \wdelta_1(L,M)
	 \fr{T_1(L,A,q)}{(q)_L}
 \tag2.28
 $$
 is now needed to iterate further.

 To find a $\left(\wgamma_1,\wdelta_1\right)$ pair we multiply equation
 (2.22) by $q^{\fr{A}{2}}\left(q^{-\fr{A}{2}}\right)$ and then replace
 $A$ by $A + 1 (A - 1)$ to get
 $$
 \align
	 \sum_{L=A+1}^M \fr{(-1)_L}{(q)_L} q^{\fr{L+A+1}{2}} 
	 T_0(L,A+1,q) & = \fr{(-1)_{M+1}}{(q)_M} \;
	 \fr{T_1(M,A+1,q)}{1 + q^{-1-A}}  \tag2.29     \\
	 \sum_{L=A-1}^M \fr{(-1)_L}{(q)_L} q^{\fr{L-A+1}{2}} 
	 T_0(L,A-1,q) & = \fr{(-1)_{M+1}}{(q)_M} \;
	 \fr{T_1(M,A-1,q)}{1 + q^{-1+A}}.  \tag2.30
 \endalign
 $$
 Adding (2.29) and (2.30) and using (2.9), (2.13) gives
 $$
 \split
	 \sum_{L=A-1}^M \fr{(-1)_L}{(q)_L} \{T_1(L+1,A,q) -
	 T_1(L,A,q)\} =   \\
	 = \fr{(-1)_{M+1}}{(q)_M}  \left\{
	 \fr{T_1(M,A+1,q)}{1 + q^{-1-A}} + \fr{T_1(M,A-1,q)}{1 + q^{-1+A}}
	 \right\}\;.
 \endsplit
 \tag2.31
 $$
 Next we treat the sum in (2.31) as follows
 $$
 \align
	 & \sum_{L=A-1}^M \fr{(-1)_L}{(q)_L} \{T_1(L+1,A,q) - T_1
	 (L,A,q)\}  =   \\
	 & = \sum_{L=A-1}^M \left\{\fr{(-1)_L}{(q)_L} - 
		 \fr{(-1)_{L+1}}{(q)_{L+1}}\right\} T_1(L+1,A,q) +
		 \tag2.32    \\
	 & + \sum_{L=A-1}^M \left\{\fr{(-1)_{L+1}}{(q)_{L+1}} 
		 T_1(L+1,A,q) - \fr{(-1)_L}{(q)_L} T_1(L,A,q)\right\} 
	 \\
	 & = - \sum_{L=A}^{M+1} \left( -q^{-1}\right)_L q^L \fr{T_1(L,A,q)}{(q)_L}
		 + \fr{(-1)_{M+1}}{(q)_{M+1}} T_1(M + 1,A,q).
 \endalign
 $$
 Combining (2.31), (2.32) and replacing $M$ by $M - 1$ yields
 $$
 \align
	 \sum_{L=A}^M \left(-q^{-1}\right)_L q^L \fr{T_1(L,A,q)}{(q)_L}
	 & = \fr{(-1)_M}{(q)_M} \bigg\{T_1(M,A,q)     \tag2.33
	 \\
	 & \qquad - \fr{(1-q^M)}{1 + q^{-1-A}} T_1(M - 1,A + 1,q)  \\
	 & \qquad -\left. \fr{(1-q^M)}{1 + q^{-1+A}} T_1(M - 1,A - 1,q)
	 \right\}  
 \endalign
 $$
 which is nothing else but (2.28) with
 $$
 \align
	 \wdelta_1(L,M) & = \theta(L \leq M) \left(- \fr1{q}\right)_L q^L
		 \tag2.34     \\
	 \wgamma_1(L,M) & = \fr{(-1)_M}{(q)_M} \bigg\{T_1(M,A,q)  
		 \tag2.35   \\ 
	 & \qquad - \fr{(1 - q^M)}{1 + q^{-1-A}} T_1(M-1,A+1,q)  \\
	 & \qquad - \left.\fr{(1 - q^M)}{1 + q^{-1+A}} T_1(M-1,A-1,q)
		 \right\}.  \\
 \endalign
 $$
 The proof of the second statement of Theorem 1 (1.13) follows again
 by the Bailey Transform argument (2.27) (with subindex $0$ replaced
 by $1$, everywhere).  Unlike (1.12), equation (1.13) does not appear
 to be a defining relation for the new Bailey pair and therefore can
 not be iterated further.

 Finally, letting $L$ tend to infinity in (1.12), (1.13) and using
 the limiting formula (2.15) gives:

 \proclaim
 {Theorem 2}  If a pair of sequences $(\walpha_{a=0,1},\wbeta_{a=0,1})$
 is subject to identities (1.11) then 
 $$
 \align
	 \sumlo (-1)_L q^{\fr{L}{2}} \wbeta_0(L) & =
	 \fr{(-1)_{\infty}(-q)_{\infty}}{(q)_{\infty}^2} \sumro
	 \fr{\walpha_0(r)}{q^{\fr{r}{2}} + q^{-\fr{r}{2}}}
	 \tag2.36	
	 \\ \\
	 \sumlo \left(-q^{-1}\right)_L q^{L} \wbeta_1(L) & =
	 \fr{(-1)_{\infty}(-q)_{\infty}}{(q)_{\infty}^2} \sumro
	 \walpha_1(r)\left\{\fr1{1 + q^{r+1}} - \fr1{1 + q^{r-1}}
	 \right\}
	 \tag2.37
 \endalign
 $$
 hold.
 \endproclaim

 \subhead
 3. \ Applications
 \endsubhead

 \subhead
 3.1 \ Preliminaries
 \endsubhead

 Recently it was shown [25] that Bailey's lemma ``connects''
 $M(p,p+1)$ models of CFT with $N = 1$ $SM(p + 1, p + 3)$ and
 $N = 2$ $SM(p + 1,1)$ models of SCFT\footnote"{$^\dagger$}"{
 Throughout this paper notations $M(p,p'),\;N = 1\;SM(p,p')$
 $N = 2$ $SM(p,p')$ stand for models of CFT and SCFT with central 
 charges
 $$
	 1 - \fr{6(p'-p)^2}{pp'} , \fr32 \left(1 - \fr{2(p-p')^2}{pp'}\right)
	 \;,\; 3\left(1 - \fr{2p'}{p}\right) \text{ respectively.}
 $$}.
 In this section we shall demonstrate that Theorem 2 leads to very 
 different relations between these models.  We begin by collecting 
 necessary definitions and formulas.

 For $A,B \in Z$ $q$-binomial coefficients $\bmatrix A \\ B
 \endbmatrix_q$ are defined as
 $$
	 \bmatrix A \\ B\endbmatrix_q = \cases
		 \fr{(q)_A}{(q)_B (q)_{A-B}} \qquad & \text{ for }
			 0 \leq B \leq A   \\
		 0  \qquad &  \text{ otherwise }.
	 \endcases
 \tag3.1
 $$
 The following properties of $q$-binomials
 $$
 \align
	 \bmatrix A \\ B \endbmatrix_{1/q} & = q^{B(B-A)}
		 \bmatrix A \\ B \endbmatrix_q    \tag3.2    \\
	 \lim_{A\rightarrow\infty} \bmatrix A \\ B  \endbmatrix_q
	 & = \fr1{(q)_B}  \tag3.3
 \endalign
 $$
 are well known.  Next we state some bosonic character formulas
 $$
	 M(p,p')[35-37]:\chi_{r,s}^{p,p'}(q) = \fr1{(q)_{\infty}}
	 \sum_{j=-\infty}^{\infty} \{q^{j(jpp' + rp'-sp)} - 
		 q^{(jp+r)(jp'+ s)}\}
 \tag3.4
 $$
 where $p' > p \geq 2$ are positive coprime integers and $r \in
 \{1,2,\dots,p-1\}$, $s \in \{1,2,\dots,p'-1\}$ are labels of
 irreducible highest weight representations
 $$
 \gathered
	 N = 1 \; SM(\wp,\wp') [37-40] : \wchi^{\wp,\wp'}_{\wr,\ws}(q)
	 = \fr{(-q^{\in_{(\wr-\ws)}})_{\infty}}{(q)_{\infty}}   \\
	 \sum_{j=-\infty}^{\infty}
	 \left\{q^{ \fr{j(j\wp\wp' + \wr\wp' - \ws \wp)}{2}} - q^{
	 \fr{(j\wp + \wr)(j\wp' + \ws)}{2}}\right\}
 \endgathered
 \tag3.5
 $$
 where
 $$
	 \in_a = \cases
		 1/2 \qquad & \text{ for } a = 0 \pmod{2}   \\
		 1 \qquad & \text{ for } a = 1 \pmod{2}   
	 \endcases
 \tag3.6
 $$
 $\wp' > \wp \geq 2$ are positive integers (with $\fr{p' - p}{2}$ and
 $p$ being coprime integers) and $\wr \in \{1,2,\dots,\wp - 1\}$,
 $\ws \in \{1,2,\dots,\wp' - 1\}$.
 $$
 \gathered
	 N = 2\;SM(\wtp,1)^{\dagger)} [41-43]: \quad
	 \wtchi_{\wtr,\wts}^{\wtp,1}(q,y) = \fr{(-q^{\wtin} y)_{\infty}
		 (-q^{\wtin}y^{-1})_{\infty}}{(q)_{\infty}^2} *
	 \\
	 \sum_{j=-\infty}^{\infty} q^{j^2\wtp + j (\wtr + \wts)}
		 \fr{1 - q^{2\wtp j + \wtr + \wts}}{(1 + y^{-1}q^{\wtp j + \wtr}			)(1 + y q^{\wtp j + \wts})}
 \endgathered
 \tag3.7
 $$
 where $\wtp \geq 2$ is a positive integer and
 \roster
 \item  in the $A$ sector, $\wtr,\wts$ are half integers with
 $0 < \wtr,\wts,\;\wtr + \wts \leq \wtp - 1$; $\wtin = 1/2$

 \item in the $P$ sector, $\wtr,\wts$ are integers with
 $0 < \wtr - 1,\wts, \; \wtr + \wts \leq \wtp - 1$; $\wtin = 1$
 \endroster
 \footnote"{}"{$^\dagger)$See [44] for the latest discussion 
 regarding (3.7).}

 All $N = 2$ $SM(\wtp,\wtp' > 1)$ characters were calculated by Ahn
 et al [45] in terms of fractional level string functions.  However,
 for the vacuum sector for $N = 2$ $SM(\wtp,\wtp > 1)$ model (with
 $\wtp > \wtp' \geq 2$; $\wtp,\wtp'$ coprime) character formula similar 
 to (3.7)
 $$
 \align
	 \wtchi_{\wtr,\wts}^{\wtp,\wtp'}(q,y) = & 
	 \fr{(-q^{1/2}y)_{\infty}(-q^{1/2}y^{-1})_{\infty}}{(q)_{\infty}^2}
	 *    \tag3.8   \\
	 \sum_{j=-\infty}^{\infty} &\fr{q^{j^2\wtp\wtp' + j(\wtr + \wts)
		 \wtp'}(1 - q^{2\wtp j+ \wtr + \wts})}{(1 + y^{-1} q^{\wtp j 
		 + \wtr})(1 + y q^{\wtp j + \wts})}\;;
	 \wtr = \wts = 1/2
 \endalign
 $$
 was recently found in [46].  Presumably, (3.8) also holds for
 sufficiently small $\wtr,\wts \in Z + \fr12$, such that the embedding
 diagram is the same as in vacuum case $\wtr = \wts = 1/2$.
 There are many important differences between $N = 2$ $SM(p,1)$ 
 and $N = 2$ $SM(p,p' > 1)$ models.  In particular, in contrast
 to the $N = 2$ $SM(p,1)$ case $N = 2$ $SM(p,p' > 1)$ models are
 neither unitary nor rational [45, 46].  Moreover, while characters
 (3.7) have nice modular properties [47], those of $N = 2$ $SM(p,p' > 1)$
 do not.  Nevertheless, one can show that characters (3.8) are, in fact,
 mock theta functions$^{\ddagger}$, i.e. they exhibit sharp asymptotic
 behaviour when $q(|q| < 1)$ tends to a rational point of unit circle.
 \footnote"{}"{$^{\ddagger}$Notion of mock theta function was introduced
 by Ramanujan in his last letter to Hardy, dated January 1920.}

 While bosonic characters (3.4, 3.5, 3.7) were known for quite some
 time, new fermionic expressions for these characters became available
 only in the last few years.  Existence of the fermionic representations
 suggests that the Hilbert space of (S)CFT can be described in terms of
 quasi-particles obeying Pauli's exclusion principle.  The equivalence
 of the bosonic and fermionic character formulas gives rise to many 
 new $q$-series identities of Rogers-Ramanujan type.  Remarkably, in
 many known cases these identities admit polynomial analogues which
 can be written as defining relations (1.11) for trinomial Bailey
 pairs.

 \subhead
 3.2 \ Polynomial analogues of generalized G\"ollnitz-Gordon identities.
 $\bold{N = 1\;SM(2,4\nu)}$, $\bold{N = 2\; SM(4\nu,1)}$ relation
 \endsubhead

 Many polynomial Fermi-Bose character identities for $N = 1$ 
 $SM(2,4\nu)$, $\nu \geq 2$ were derived in [33, 34].  Not to
 overburden our narrative with cumbersome notations we shall 
 consider here only the simplest of these identities
 $$
 \gathered
	 \sum_{n_1,\dots,n_{\nu} = 0}^{\infty} q^{\fr{n_1^2}{2}- n_1 N_2
	 + \overset{\nu}\to{\underset{j=2}\to\sum}  
		 N_j^2} \bmatrix N_2 \\ n_1 \endbmatrix_q
	 \prod_{i=2}^{\nu} \bmatrix n_i + L + n_1 - 2 
	 \overset{i}\to{\underset{j=2}\to\sum}  N_j \\ n_i \endbmatrix =
	 \\
	 = \sum_{j=-\infty}^{\infty} (- \ )^j q^{\nu j^2 + j/2} \{T_0
	 (L,2\nu j,q) + T_0(L,2\nu j + 1,q)\}
 \endgathered
 \tag3.9
 $$
 where
 $$
	 N_j = n_j + n_{j+1} + \dots + n_{\nu}
 \tag3.10	
 $$
 Letting $L$ in (3.9) tend to infinity yields
 $$
 \gathered
	 \sum_{n_1,\dots,n_{\nu}=0}^{\infty} q^{\fr{n_1^2}{2}- n_1 N_2
	 + \overset{\nu}\to{\underset{j=2}\to\sum}  N_j^2}  
		 \bmatrix N_2 \\ n_1 \endbmatrix_q
	 \fr1{(q)_{n_2} \cdots (q)_{n_{\nu}}} =  \\
	 = \fr{(-q^{1/2})_{\infty}}{(q)_{\infty}} \sum_{j=-\infty}^{\infty}
	 (- \ )^j q^{\nu j^2 + j/2} = \wchi_{1,2\nu - 1}^{2,4\nu}(q),
 \endgathered
 \tag3.11
 $$
 where we used (3.3) and a limiting formula
 $$
	 \lim_{L\rightarrow\infty} \{T_0(L,A,q) + T_0(L,A + 1,q)\}
	 = \fr{(-q^{1/2})_{\infty}}{(q)_{\infty}}
 \tag3.12
 $$
 proven in [28].  Identity (3.11) is nothing else but Andrews
 generalization of G\"ollnitz-Gordon identities [48].  A moment's 
 reflection shows that (3.9) is in the form (1.11) with
 $$
	 \walpha_0(r) = \cases
		 (-1)^j q^{\nu j^2}(q^{j/2} + q^{-j/2}) \qquad & 
			 \text{ for } r = 2 \nu j,\quad j > 0  \\
		 1 \qquad & \text{ for } r = 0   \\
		 (-1)^j q^{\nu j^2 + j/2} \qquad & \text{ for }
			 r = 2\nu j + 1,\quad j \geq 0   \\
		 (-1)^j q^{\nu j^2 - j/2} \qquad & \text{ for }
			 r = 2\nu j - 1,\quad j \geq 1
	 \endcases
 \tag3.13
 $$ 
 $$
	 \wbeta_0(L) = \fr1{(q)_L} \sum_{n_1,\dots,n_{\nu}=0}^{\infty}
	 q^{\fr{n_1^2}{2} - n_1 N_2 + \overset{\nu}\to{\underset{j=2}\to\sum}
		  N_j^2}
	 \bmatrix N_2 \\ n_1 \endbmatrix_q \prod_{j=2}^{\nu}
	 \bmatrix n_i + L + n_1 - 2 \overset{i}\to{\underset{j=2}\to\sum} 
		 N_j  \\ n_i\endbmatrix_q.
 \tag3.14
 $$
 Substituting (3.13), (3.14) into (2.36) gives
 $$
 \gathered
	 \sum_{L,n_1,\dots,n_{\nu}=0}^{\infty} \fr{(-1)_L}{(q)_L} 
	 q^{\fr{L+n_1^2}{2} - n_1 N_2 + \overset{\nu}\to{\underset{j=2}\to\sum}
		   N_j^2}
	 \bmatrix N_2 \\ n_1\endbmatrix_q  \prod_{i=2}^{\nu}
	 \bmatrix n_i + L + n_1 - 2\overset{i}\to{\underset{j=2}\to\sum} 
		  N_j  \\ n_i\endbmatrix_q	=
	 \\
	 = \fr{(-1)_{\infty}(-q)_{\infty}}{(q)_{\infty}^2} 
		 \sum_{j=-\infty}^{\infty} (- \ )^j q^{\nu j^2 + j/2}
		 \left\{\fr{q^{\nu j}}{1 + q^{2\nu j}} +
		 \fr{q^{\nu j + 1/2}}{1 + q^{2\nu j + 1}}\right\} =
	 \\
	 = \wtchi_{1/2,2\nu + 1/2}^{4\nu,1}\; (q,q^{1/2}) + q^{1/2}
	 \wtchi_{3/2,2\nu - 1/2}^{4\nu,1} (q,q^{1/2})
 \endgathered
 \tag3.15
 $$
 which establishes advertised relation between $N = 1$ $SM(2,4\nu)$ and
 $N = 2$ $SM(4\nu,1)$ models of SCFT.  Moreover, left hand side of
 equation (3.15) provides new fermionic companion form for $N = 2$
 $SM(4\nu,1)$ characters.  This form is quite different from known
 fermionic representation [21,25] given in terms of $\D_{4\nu}$-Cartan
 matrix.

 \subhead
 3.3. \ Trinomial Bailey flow from $\bold{M(3,4)}$ (Ising) model to 
 $\bold{N = 2}$ $\bold{SM(6,1)}$ model of SCFT
 \endsubhead

 In [29] the following polynomial identity
 $$
 \align
	 \sum_{j=0}^L \bmatrix L \\ j \endbmatrix_q q^{j^2/2}
	 & = \sum_{j=-\infty}^{\infty} q^{6j^2 + j}(T_0(L,6j,q) +
	 T_0(L,6j + 1,q))  \tag3.16
	 \\
	 & - \sum_{j=-\infty}^{\infty} q^{6j^2 + 5j+1}(T_0(L,6j+2,q) +
	 T_0(L,6j + 3,q))  
 \endalign
 $$
 was proven.  One may check that in the limit $L \to \infty$ this
 identity reduces to Fermi-Bose character identity for $M(3,4)$ 
 (Ising) model
 $$
 \align
	 \sum_{j\geq 0}  \fr{q^{j^2/2}}{(q)_j} & =  \fr{(-q^{1/2})_{\infty}}
		 {(q)_{\infty}} \sum_{j=-\infty}^{\infty} (q^{6j^2 + j}-
		 q^{6j^2 + 5j+1})   \\
	 & = \chi_{1,1}^{3,4}(q) + q^{1/2}\; \chi_{2,1}^{3,4}(q).
 \tag3.17
 \endalign
 $$
 The middle expression in (3.17) is remarkably  similar to (3.5) with $\wp=3,\;\wp'=4$.
 This similarity suggests an interpretation of (3.17) as a character of some extended 
 Virasoro algebra. 

 It is straightforward to verify that (3.16) is the defining relation
 (1.11) for trinomial pair 
 $$
	 \walpha_0(r) = \cases 
		 q^{6j^2}(q^j + q^{-j}) \qquad & \text{ for } r = 6j,
			 \quad j > 0  \\
		 1 \qquad & \text{ for } r = 0   \\
		 q^{6j^2 +j} \qquad & \text{ for } r = 6j + 1, j \geq 0  \\
		 q^{6j^2 -j} \qquad & \text{ for } r = 6j - 1, j > 0  \\
		 - q^{6j^2 + 5j + 1} \qquad & \text{ for } r = 6j + 2
			 \text{ and } r = 6j + 3\,, \quad j \geq 0   \\
		 - q^{6j^2 - 5j + 1} \qquad & \text{ for } r = 6j - 2
			 \text{ and } r = 6j - 3\,, \quad j > 0   
	 \endcases 
 \tag3.18
 $$
 $$
	 \wbeta_0(L) = \fr1{(q)_L}   \sum_{j \geq 0} q^{\fr{j^2}{2}}
	 \bmatrix{L}  \\  {j}\endbmatrix_q.
 \tag3.19
 $$
 Next we apply Theorem 2 to the pair (3.18, 3.19), the result is
 $$
 \align
	 \sum_{L,j\geq 0} q^{\fr{L + j^2}{2}} \fr{(-1)_L}{(q)_L}
	 \bmatrix L \\ j \endbmatrix_q & = \fr{(-1)_{\infty}(-q)_{\infty}}
		 {(q)_{\infty}^2} \left\{ \sum_{j=-\infty}^{\infty}
	 q^{6j^2 + j} \left( \fr{q^{3j}}{1 + q^{6j}} + 
		 \fr{q^{3j + 1/2}}{1 + q^{6j+1}}\right)  \right.  
	 \\
	 & \hskip .5in \left. -\sum_{j=-\infty}^{\infty} 
		 q^{6j^2 + 5j + 1} \left( \fr{q^{3j + 1}}{1 + q^{6j + 2}} + 
		 \fr{q^{3j + 3/2}}{1 + q^{6j+3}}\right)\right\}
	 \\
	 & =\wtchi_{\fr12,\fr32}^{6,1} \left(q,q^{\fr12}\right) +
		 q^{\fr12} \wtchi_{\fr12,\fr32}^{6,1} \left(q,q^{\fr32}\right).
		 \tag3.20
 \endalign
 $$
 Recently, Warnaar proposed polynomial identities similar to (3.16) 
 for all models $M(p,p + 1)$, $p \geq 3$ [16].  $p = 3$ Case is the
 one treated above.  We have checked that his conjecture implies 
 the following identities
 $$
 \align
	 \left\{  \sum_{L,m_1,\dots,m_{p-2}\geq 0}\right.  &   q^{\fr{L + m_1^2}{2} - Lm_1
	 + \fr{m_1 m_2}{4} + \fr14 \overset{p-2}\to{\underset{i=2}\to\sum}
		  m_i(2m_i - m_{i-1}
	 - m_{i+1})}     \\
	     & \left.  \fr{(-1)_L}{(q)_L} \bmatrix L \\ m_1\endbmatrix_q
	  \prod_{i=2}^{p-2} 
		\bmatrix \fr{m_{i-1} + m_{i+1}}{2}  \\ m_i \endbmatrix_q
		\right\}	=  \tag3.21
	\\
	& = q \fr{(-1)_{\infty}(q)_{\infty}}{(q)_{\infty}^2}
		\sum_{j=-\infty}^{\infty}
	q^{j^2 p(p-1) + j(p-1)}  \fr{1 - q^{4pj+2}}{(1 + q^{2pj})
		(1 + q^{2pj+2})}
\endalign
$$
where $m_{p-1}\equiv 0$.  Therefore Trinomial Bailey flow for 
$p = 1 \pmod{2}$ is
$$
	M(p, p+1) \longleftrightarrow N = 2 \; SM\left(2p,
		\fr{p-1}{2}\right).
\tag3.22
$$
This is to be contrasted with Bailey flow discussed in [25] 
where one has
$$
	M(p, p+1) \longleftrightarrow N = 2\; SM(p + 1,1).
\tag3.23
$$

\subhead
3.4 \ Results related to Rogers-Ramanujan identities
\endsubhead

It is well known that Rogers-Ramanujan identities
$$
\aligned
	\sum_{j\geq 0} \fr{q^{j(j+a)}}{(q)_j} & = \fr1{(q)_{\infty}}
	\sum_{j=-\infty}^{\infty}  \left\{ q^{j(10j + 1 + 2a)} 
	- q^{(2j+1)(5j +2 - a)}\right\}   \\
	& = \chi_{1,2-a}^{2,5} (q)\;;\; a = 0,1  
\endaligned
\tag3.24
$$
admit polynomial analogues
$$
\aligned
	\sum_{j \geq 0} q^{j(j+a)} \bmatrix 2L - j - a  \\ j \endbmatrix_q
	= \sum_{j=-\infty}^{\infty} & \left\{ q^{j(10j + 1 + 2a)}
	\bmatrix 2L  \\  L - 5j - a\endbmatrix_q  \right.  \\
	& - \left. q^{(2j + 1)(5j + 2-a)} \bmatrix 2L \\ L - 5j - 2
		\endbmatrix_q	\right\}
\endaligned
\tag3.25
$$
which reduce to (3.24) as $L \to \infty$.  It is rather surprising
that polynomials appearing in (3.25) have $q$-trinomial 
representation as well [29].  In particular, for $a = 0$, one
has
$$
\aligned
	\sum_{j\geq 0} q^{j2} \bmatrix 2L - j  \\ j\endbmatrix_q
	= \sum_{j=-\infty}^{\infty} & \left\{ q^{60j^2 - 4j}
	\pmatrix L,10j;q^2  \\ 10j \endpmatrix_2  \right.  \\
	& - q^{60j^2 + 44j +8} \pmatrix L,10j + 4;q^2  \\ 10j + 
		4\endpmatrix_2     \\
	& + q^{60j^2 + 16j + 1} \pmatrix L,10j + 1;q^2  \\ 10j + 
		1\endpmatrix_2     \\
	& - \left. q^{60j^2 + 64j + 17} \pmatrix L,10j+5; q^2  \\
		10j + 5 \endpmatrix_2\right\}
\endaligned
\tag3.26
$$
Identity (3.26) is not of the form (1.11).  However, if we replace
$q$ by $\fr1{\sqrt{q}}$ in (3.26) and multiply the result by
$q^{L^2/2}$ we obtain with the help of (2.6), (3.2)
$$
\aligned
	\sum_{j\geq 0} q^{\fr{j^2}{2}} 
		\bmatrix L + j \\ 2j \endbmatrix_{\sqrt{q}}
	& = \sum_{j=-\infty}^{\infty} q^{20j^2 + 2j} (T_0(L,10j,q)
		+ T_0(L,10j+1,q))
	\\
	& - \sum_{j=-\infty}^{\infty} q^{20j^2 + 18j + 4}(T_0(L,10j 
		+ 4,q) + T_0(L,10j+5,q))
\endaligned
\tag3.27
$$
which gives rise to trinomial Bailey pair
$$
	\wbeta_0(L) = \fr1{(q)_L} \sum_{j\geq 0}q^{\fr{j^2}{2}}
	\bmatrix L + j \\ 2j \endbmatrix_{\sqrt{q}}
\tag3.28
$$
\vskip -.2in
\baselineskip 12pt
$$
	\walpha_0(r) =  \left\{
\aligned
		q^{20j^2}(q^{2j} + q^{-2j})  \qquad & \text{ for }
			r = 10j,j > 0   \\
		1 \qquad & \text{ for } r = 0   \\	
		q^{20j^2 + 2j} \qquad & \text{ for } r = 10j + 1,
			j \geq 0     \\
		q^{20j^2 - 2j} \qquad & \text{ for } r = 10j - 1,
			j > 0    \\
		-q^{20j^2 + 18j + 4} \qquad & \text{ for } r = 10j
			+ 4 \text{ and } r = 10j + 5\,, \quad j \geq 0
		\\	
	-q^{20j^2 - 18j + 4} \qquad & \text{ for } r = 10j
			- 4 \text{ and } r = 10j - 5\,, \quad j > 0
\endaligned  \right.	
 \tag3.29   
$$

\baselineskip 20pt
Next we apply Theorem 2 to derive
$$
\multline
	\sum_{L,j\geq 0} q^{\fr{L + j^2}{2}} \fr{(-1)_L}{(q)_L}
	\bmatrix L + j  \\ 2j \endbmatrix_{\sqrt{q}}   =     \\
	= \fr{(-1)_{\infty}(-q)_{\infty}}{(q)_{\infty}^2} \left\{
		\sum_{j=-\infty}^{\infty} q^{20j^2 + 2j} \left(
		\fr{q^{5j}}{1 + q^{10j}} + \fr{q^{5j + \fr12}}
			{1 + q^{10j+1}} \right)\right.        \\
	- \left.\sum_{j=-\infty}^{\infty} q^{20j^2 + 18j + 4}\left(
	\fr{q^{5j + 2}}{1 + q^{10j + 4}} + \fr{q^{5j + \fr52}}
			{1 + q^{10j+5}} \right)\right\}.        
\endmultline
\tag3.30
$$
We note that the expression on the right hand side of (3.30) bears a
strong resemblance to formula (3.8) with $\wp=10,\;\wp'=2$. It is 
also similar in form to the $\Phi(q)$ and $\Psi(q)$ considered by
Ramanujan in his development of the fifth-order mock theta functions [49].

Expressions in (3.30) are not modular functions. Nevertheless, using
Poisson summation formula one can show that asymptotic behaviour of (3.30)
can be neatly expressed in terms of exponential forms.  For instance,
when $q = e^{-t}$ and $t \to 0$ we have for (3.30)
$$
	\sqrt{\fr{2}{5 \pi t}}  \cos \fr{\pi}{10} 
	e^{\fr{9 \pi^2}{20 t}}
\tag3.31
$$
Proof of (3.31) along with asymptotic analysis at nonunitary 
characters (3.8) will be given elsewhere.

\subhead
4. \ Discussion
\endsubhead

It is widely believed that different fermionic expressions for the
(super) conformal character are related to different integrable
pertrubations of the same (super) conformal model.  Thus, it would
be interesting to identify perturbations which correspond to the
new fermionic representations for $N = 2$ $SM(4\nu,1)$ characters
found in Section 3.2.

Furthermore, following [25, 26], it is tempting to interpret
the relations
$$
\gathered
	N = 1\; SM(2,4\nu) \longleftrightarrow N = 2\; SM(4\nu,1)
	\\
	M(p,p + 1) \longleftrightarrow N = 2 \; SM\left(2p, 
		\fr{p - 1}{2}\right)\;, \; p = 1 \pmod{2}
\endgathered
\tag4.1
$$
established here as massless renormalization group flows.  If such
an interpretation is indeed correct, then one should be able to 
carry out Thermodynamic Bethe Ansatz (TBA) analysis of these flows along
the lines of [50]. We expect that related TBA systems will have the same
incidence structure as that of fermionic forms discussed in section 3.
Also, we would like to point out that a ``folding
in half'' relation between $N = 2$ $SM(4\nu,1)$ and $N = 1$ $SM(2,4\nu)$
models has already been noticed in [51,52]. Partition theoretical 
interpretation of our results will undoubtably lead to construction of
subtractionless bases for $N=2$ super Virasoro modules. 

From the mathematical point of view it is highly desirable to find
an appropriate $q$-hypergeometric background for the Trinomial analogue
of Bailey's lemma.  Recall that the classical Bailey's lemma is
intimately related to the $q$-Pfaff-Saalsch\"utz formula
$$
	\sum_{j=0}^n  \;\fr{(q^{-n})_j (c/a)_j(c/b)_j}{(q)_j (c)_j
		(cq^{1-n}/ab)_j}\;q^j = \fr{(a)_n(b)_n}{(c)_n(ab/c)_n}
\tag4.2
$$
which was first derived by Jackson [53]. In this direction we have already
determined that
$$
	\sumlo (-q^{-n})_L q^{\fr{(1+2n)L}{2}} \wbeta_0(L) =
	\fr{(-q^{n+1})_{\infty}^2}{(q)_{\infty}(q^{2n+1})_{\infty}} \sumro
	\fr{q^{\fr{1+2n}{2}r}(-q^{-n})_r}{(-q^{n+1})_r}\walpha_0(r)
	\tag4.3
$$
with $n=0,1,2,3,...$ . Identity (4.3) can be derived from (2.12, 2.13, 2.15, 2.22) after a bit of labour.
It immediately follows that $q^n$ in (4.3) may be replaced by an arbitrary parameter,
say $\rho$. Details will be given elsewhere [54].

Building on a proposal made in [31], Schilling and Warnaar defined
and extensively studied $q$-multinomials [55-57].  One may wonder
if these new objects will lead to additional generalizations of 
Bailey's lemma.  We strongly believe that the answer is ``yes''
and hope to say more about it in a subsequent paper.

\subhead
Note Added
\endsubhead

Soon after this paper was completed, Warnaar provided a simple and elegant
proof of the conjecture from [16] used in deriving (3.21). Moreover,
he has shown that each ordinary Bailey pair gives rise to a trinomial
Bailey pair.  In particular, he demonstrated that the trinomial Bailey pair
(3.18-19) is a ``descendant'' of the A(1) and A(2) Bailey pairs of Slater's
list [4].

\subhead
Acknowledgement
\endsubhead

We would like to thank B.M. McCoy, Z. Reti, K. Voss, and S.O. Warnaar for
interesting discussions and helpful comments.
This work was partially supported by NSF grant: DMS-9501101.


\Refs

\ref
  \no 1
  \by W. N. Bailey
  \paper
  \jour Proc. London Math. Soc. (2)
  \vol 49
  \yr 1947
  \pages 421	
\endref

\ref
  \no 2
  \by W. W. Bailey
  \paper 
  \jour Proc. London Math. Soc (2)
  \vol 50
  \yr 1949
  \pages 1	
\endref

\ref
  \no 3
  \by F. J. Dyson 
  \paperinfo in Ramanujan Revisited ed. by G. E. Andrews et al
	(Academic Press, 1988) p. 7
\endref

\ref
  \no 4
  \by L. J. Slater
  \paper
  \jour Proc. London Math. Soc. (2)
  \vol 53
  \yr 1951
  \pages 460	
\endref

\ref
  \no 5
  \by L. J. Slater
  \paper
  \jour Proc. London Math. Soc. (2)
  \vol 54
  \yr 1952
  \pages 147	
\endref

\ref
  \no 6
  \by G. E. Andrews
  \paper 
  \jour Pac. Journ. Math.
  \vol 114
  \yr 1984
  \pages 267	
\endref

\ref
  \no 7
  \by A. K. Agarwal, G. E. Andrews, D. M. Bressoud
  \paper
  \jour Jour. Indian Math. Soc.
  \vol 51
  \yr 1987
  \pages 57
\endref

\ref
  \no 8
  \by D. M. Bressoud
  \paperinfo in Ramanujan Revisited ed. by G. E. Andrews et al
		(Academic Press, 1988) p. 57
\endref

\ref
  \no 9
  \by G. E. Andrews
  \paper ``$q$-series:  Their development and application in
	analysis, number theory, combinatorics, physics and computer
	algebra''
  \paperinfo (American Math. Society, Providence, Rhode Island,
		1986)
\endref

\ref
  \no 10
  \by S. C. Milne, G. M. Lilly
  \paper
  \jour Bull. Amer. Math. Soc.
  \vol 26
  \yr 1992
  \pages 258	
\endref

\ref
  \no 11
  \by E. Melzer
  \paper
  \jour Int. Jour. Mod. Phys.
  \vol A9
  \yr 1994
  \pages 1115
\endref

\ref
  \no 12
  \by A. Berkovich
  \paper
  \jour Nucl. Phys.
  \vol B431
  \yr 1994
  \pages 315	
\endref

\ref
  \no 13
  \by O. Foda, Y.-H. Quano
  \paper
  \jour Int. Jour. Mod. Phys.
  \vol A10
  \yr 1995
  \pages 2291	
\endref

\ref
  \no 14
  \by A. N. Kirillov
  \paper
  \jour Prog. Theor. Phys. Suppl.
  \vol 118
  \yr 1995
  \pages 61	
\endref

\ref
  \no 15
  \by S. O. Warnaar
  \paper
  \jour Jour. Stat. Phys.
  \vol 82
  \yr 1996
  \pages 657
\endref

\ref
  \no 16
  \by S. O. Warnaar
  \paper
  \jour Jour. Stat. Phys.
  \vol 84
  \yr 1996 
  \pages 49	
\endref

\ref
  \no 17
  \by A. Berkovich, B. M. McCoy
  \paper
  \jour Lett. Math. Phys.
  \vol 37
  \yr 1996
  \pages 49
\endref

\ref
  \no 18
  \by A. Berkovich, B. M. McCoy, A. Schilling
  \paper ``Rogers-Schur-Ramanujan type identities for $M(p,p')$
	minimal models of conformal field theory
  \paperinfo $q$-alg/9607020, submitted to Commun. Math. Phys.
\endref

\ref
  \no 19
  \by R. Kedem, T. R. Klassen, B. M. McCoy, E. Melzer
  \paper
  \jour Phys. Lett. 
  \vol B304
  \yr 1993	
  \pages 263	
\endref

\ref
  \no 20
  \by R. Kedem, T. R. Klassen, B. M. McCoy, E. Melzer
  \paper
  \jour Phys. Lett. 
  \vol B307
  \yr 1993
  \pages 68	
\endref

\ref
  \no 21
  \by S. Dasmahapatra, R. Kedem, T  R. Klassen, B. M. McCoy,
	E. Melzer
  \paper
  \jour Int. Jour. Mod. Phys.
  \vol B7
  \yr 1993
  \pages 3677	
\endref

\ref
  \no 22
  \by O. Foda, Y.-H. Quano
  \paper ``Virasozo character identities from Andrews-Bailey
	construction''
  \jour hep - th/9408086, Int. Jour. of Mod. Phys. A 
  \toappear
\endref

\ref
  \no 23
  \by A. Schilling, S. O. Warnaar
  \paper
  \jour Int. Jour. Mod. Phys.
  \vol B11
  \yr 1997
  \pages 189	
\endref

\ref
  \no 24
  \by A. Schilling, S. O. Warnaar
  \paper ``A higher level Bailey lemma: Proof and Applications''
  \paperinfo $q$-alg/9607019, Ramanujan Jour. (to appear)
\endref

\ref
  \no 25
  \by A. Berkovich, B. M. McCoy, A. Schilling
  \paper
  \jour Physica
  \vol A228
  \yr 1996
  \pages 33	
\endref

\ref
  \no 26
  \by L. Chin
  \paper ``Central charge and the Andrews-Bailey Constructions''
  \paperinfo hep-th/9607168
\endref

\ref 
  \no 27
  \by A. Berkovich, B. M. McCoy, A. Schilling, S. O. Warnaar
  \paper ``Bailey flows and Bose-Fermi identities for the
   conformal coset models
   $(A_1^{(1)})_N\times(A_1^{(1)})_{N'}/(A_1^{(1)})_{N+N'}$''
  \paperinfo hep-th/9702026
\endref


\ref
  \no 28
  \by G. E. Andrews, R. J. Baxter 
  \paper
  \jour Jour. Stat. Phys. 
  \vol 47
  \yr 1987
  \pages 297	 
\endref

\ref
  \no 29
  \by G. E. Andrews
  \paper
  \jour Jour. Amer. Math. Soc.
  \vol 3
  \yr 1990
  \pages 653	
\endref

\ref
  \no 30
  \by G. E. Andrews
  \paperinfo in Analytic Number Theory, B. Berndt et al eds.,
	Birkh\"auser, Boston, 1990, pp. 1--11
\endref

\ref
  \no 31
  \by G. E. Andrews
  \paper 
  \jour Contemp. Math.
  \vol 166
  \yr 1994
  \pages 141
\endref

\ref
  \no 32
  \by S. O. Warnaar, P. A. Pearce
  \paper
  \jour Jour. Phys.
  \vol A27
  \yr 1994
  \pages L891	
\endref

\ref
  \no 33
  \by A. Berkovich, B. M. McCoy, W. P. Orrick
  \paper
  \jour Jour. Stat. Phys.
  \vol 83
  \yr 1996
  \pages 795 
\endref

\ref 
  \no 34
  \by A. Berkovich, B. M. McCoy
  \paper ``Generalizations of Andrews-Bressoud Identities for $N=1$
	Superconformal Model $SM(2,4\nu)$''
  \paperinfo hep-th/9508110, to appear in Int. Jour. of Math. and
	Comp. Modelling
\endref

\ref
  \no 35
  \by B. L. Feigin, D. B. Fuchs
  \paper
  \jour Funct. Anal. Appl.
  \vol 17
  \yr 1983
  \pages 241
\endref

\ref
  \no 36
  \by A. Rocha-Caridi
  \paperinfo in Vertex Operators in Mathematics and Physics, ed.
	J. Lepowsky et al (Springer, Berlin, 1985) 
\endref

\ref
  \no 37
  \by V. K. Dobrev
  \paper
  \jour Suppl. Rendiconti Circolo Matematici di Palermo, Serie II, Numero 14,
  \yr 1987
  \pages 25
\endref

\ref
  \no 38
  \by B. L. Feigin, D. B. Fuchs
  \paper
  \jour Func. Anal. Appl.
  \vol 16
  \yr 1982
  \pages 114
\endref

\ref
  \no 39
  \by P. Goddard, A. Kent, D. Olive
  \paper
  \jour Comm. Math. Phys.
  \vol 103
  \yr 1986
  \pages 105
\endref

\ref
  \no 40 
  \by A. Meurman, A. Rocha-Caridi
  \paper 
  \jour Comm. Math. Phys.
  \vol 107
  \yr 1986
  \pages 263
\endref

\ref
  \no 41
  \by V. K. Dobrev
  \paper
  \jour Phys. Lett.
  \vol B186
  \yr 1987
  \pages 43
\endref

\ref
  \no 42
  \by Y. Matsuo
  \paper 
  \jour Prog. Theor. Phys.
  \vol 77
  \yr 1987
  \pages 793
\endref

\ref
  \no 43
  \by E. B. Kiritsis
  \paper 
  \jour Int. Jour. Mod. Phys.
  \vol A3
  \yr 1988
  \pages 1871
\endref

\ref
  \no 44
  \by M. D\"orrzapf
  \paper 
  \jour Comm. Math. Phys.
  \vol 180
  \yr 1996
  \pages 195
\endref

\ref
  \no 45
  \by C. Ahn, S. Chung, S.-H.Tye
  \paper
  \jour Nucl. Phys.
  \vol B365
  \yr 1991
  \pages 191
\endref

\ref
  \no 46
  \by W. Eholzer, M. R. Gaberdiel
  \paper ``Unitarity of rational $N = 2$ superconformal theories''
  \paperinfo hep-th/9601163
\endref

\ref
  \no 47
  \by F. Ravanini, S.-K. Yang 
  \paper
  \jour Phys. Lett.
  \vol B195
  \yr 1987
  \pages 202	
\endref

\ref
  \no 48
  \by G. E. Andrews
  \paper The theory of partitions
  \paperinfo (Addison-Wesley, London, 1967)
\endref

\ref
  \no 49
  \by G. E. Andrews, F. G. Garvan
  \paper
  \jour Adv. in Math.
  \vol 73
  \yr 1989
  \pages 242
\endref

\ref
  \no 50
  \by Al. Zamolodchikov
  \paper
  \jour Nucl. Phys.
  \vol B358
  \yr 1991
  \pages 524
\endref

\ref
  \no 51
  \by E. Melzer
  \paper ``Supersymmetric Analogs of the Gordon-Andrews Identities, 
	and related TBA systems''
  \paperinfo hep-th/9412154
\endref

\ref
  \no 52
  \by M. Moriconi, K. Schoutens
  \paper 
  \jour Nucl. Phys.
  \vol B464
  \yr 1996
  \pages 472	
\endref

\ref
  \no 53
  \by F. H. Jackson
  \paper
  \jour Messenger of Math.
  \vol 39
  \yr 1910
  \pages 745
\endref

\ref
  \no 54
  \by G. E. Andrews, A. Berkovich
  \paper
  \paperinfo in preparation
\endref

\ref
  \no 55
  \by A. Schilling
  \paper 
  \jour Nucl. Phys.
  \vol B467
  \yr 1996
  \pages 247	
\endref

\ref
  \no 56
  \by S.O. Warnaar
  \paper ``The Andrews-Gordon identities and $q$-multinomial coefficients'',
	$q$-alg/9601012
  \jour Comm. Math. Phys.
  \toappear
\endref

\ref
  \no 57
  \by A. Schilling, S. O. Warnaar
  \paper ``Supernomial coefficients, polynomial identities and q-series''
  \paperinfo q-alg/9701007       
\endref

\ref
  \no 58
  \by S. O. Warnaar
  \paper ``A note on the trinomial analogue of Bailey's lemma''
  \paperinfo q-alg/9702021
\endref

\endRefs
\enddocument